\documentclass{article}

\newtheorem{theorem}{Theorem}
\newtheorem{corollary}{Corollary}

\newtheorem{proposition}{Proposition}

\newtheorem{definition}{Definition}
\newtheorem{remark}{Remark}
\newcommand{\ep}{\varepsilon}

\newcommand{\ds}{\displaystyle}
\newcommand{\ph}{\varphi}
\newcommand{\R}{{\rm I}\!{\rm  R}}
\def\R{{\rm I}\!{\rm  R}}
\def\la1{\lambda_1}
\def\grad{\nabla}
\def\ph1{\varphi_1}
\newcommand{\N}{{\bf N}}
\newcommand{\Ss}{{\bf S}_N}

\title{The Dirichlet problem for singular fully  nonlinear operators }
\date{}

\author{Isabeau Birindelli, Fran\c coise Demenge}

\begin{document}
\maketitle

{\footnotesize \centerline{}
\centerline{Universit\`a di Roma "La Sapienza"}
 \centerline{Piazzale Aldo Moro 5}
 \centerline{00185 Roma, Italy}
 \centerline{and}
 \centerline{Universit\'e de Cergy-Pontoise}
 \centerline{2 Avenue Adolphe Chauvin}
 \centerline{95302 Cergy-Pontoise, France}
}

\section{Introduction}
In the last decade, fully-nonlinear elliptic equations have been extensively studied, both in the variational and the non variational setting. The concept of viscosity solution is particularly appropriate when considering non variational fully-non linear operators. It is interesting to recall that this concept relies deeply  on a local comparison test and it is, hence, the perfect tool for studying comparison and maximum principles. 

In a different contest of linear operators and strong 
solutions, Berestycki, Nirenberg and Varadhan in \cite{BNV}, used the
maximum principle to define the concept of eigenvalue and to prove
existence of solutions for  Dirichlet problems when the boundary of the 
domain doesn't satisfy any regularity condition.

These considerations have lead to extend 
the concept of eigenvalue to the larger class of fully-nonlinear
operators, using a definition which is analogous to the one given in
\cite{BNV} via the Maximum Principle, together with the concept of 
viscosity solutions.

We shall now make this more precise. Let $\Omega$ be a bounded 
domain of $\R^N$ with a ${\mathcal C}^2$ boundary, let  $S$ be the set of
symmetric matrices 
$N\times N$, let
$\alpha>-1$,  we consider operators
$$ F(x,\grad u, D^2u)+b(x)\cdot\grad u|\grad u|^\alpha+ c(x)|u|^\alpha u$$
with
$F:\Omega\times\R^N\times S$ continuous, satisfying
\begin{itemize}
 \item[(H1)] $F(x,tp,\mu X)=|t|^{\alpha}\mu F(x,p,X)$, $\forall t\in
\R^\star$,
$\mu\in\R^+$ $\forall x\in\Omega$.

\item[(H2)] There exist $A\geq a>0$ such that for all $x\in\Omega$, $p\neq 0$ and for all $(M,N)\in S^2$, $N\geq 0$
$$a |p|^\alpha {\rm tr} N\leq F(x,p,M+N)-F(x,p,M)\leq A |p|^\alpha{\rm tr} N$$
\end{itemize}
and with $b$ and $c$ bounded and continuous.

The role of eigenvalue will be played by
$$\overline{\lambda} = \sup\{ \lambda, \exists \varphi>0\ {\rm in}\ 
\Omega,
F(x,\nabla \varphi,D^2\varphi)+b(x).\nabla \varphi|\nabla \varphi|^\alpha+
(c(x)+\lambda)|\varphi|^\alpha\varphi\leq 0\ {\rm in}\ 
\Omega\}$$
and 
$$\underline{\lambda} = \sup\{ \lambda, \exists \varphi<0\ {\rm in}\ 
\Omega , F(x,\nabla \varphi, D^2\varphi)+b(x).\nabla \varphi|\nabla
\varphi|^\alpha+ (c(x)+\lambda)|\varphi|^\alpha\varphi\geq 0\  {\rm in}\ 
\Omega\}$$

Precisely it is possible to prove that there exist two non trivial
functions $\phi_1\geq 0$  and $\phi_2\leq 0$ respectively solution of
$$
\left\{\begin{array}{lc}
F(x,\nabla \phi_1, D^2 \phi_1) +b(x)\cdot\nabla  \phi_1|\nabla \phi_1|^{\alpha}
+(c(x)+\overline\lambda)\phi_1^{\alpha+1}=0 &  \mbox{in } \  \Omega\\ 
\phi_1=0 & \mbox{on }
\partial\Omega.
\end{array}
\right.
$$
and of
$$
\left\{\begin{array}{lc}
F(x,\nabla \phi_2, D^2 \phi_2) +b(x)\cdot\nabla \phi_2|\nabla \phi_2|^{\alpha}
+(c(x)+\overline\lambda)\phi_2|\phi_2|^{\alpha}=0 &  \mbox{in } \  \Omega\\ 
\phi_2=0 & \mbox{on }
\partial\Omega.
\end{array}
\right.
$$
In the case $\alpha=0$ these results are due to Busca, Esteban, Quaas \cite{BEQ}, \cite{Q}, and
Ishii,  Yoshimura \cite{IY}, Quaas, Sirakov \cite{QS}. We wish to mention the recent work on multiplicity of solutions due to Sirakov \cite{S} for related operators and the pioneering work of P.L. Lions \cite{PLL}.

When $\alpha\neq 0$ these problems have been 
studied in 
\cite{BD2}, \cite{BD7}. In particular, for 
$\lambda<\overline\lambda$, we proved that the maximum principle holds
and that there exists a  solution $u$ for  the Dirichlet problem when the
data
$f$ is negative in $\Omega$ and for zero  boundary condition :

$$
\left\{\begin{array}{lc}
F(x,\nabla u, D^2 u) +b(x)\cdot\nabla u|\nabla u|^{\alpha}
+(c(x)+\lambda)u|u|^\alpha= f(x)&  \mbox{in } \  \Omega\\ 
u=0 & \mbox{on }
\partial\Omega.
\end{array}
\right.
$$
The scope of the present work is to enlarge these results to operators which are not homogeneous,
 to data that may change sign and non zero boundary condition.

Precisely,  for $\lambda<\min\{\overline\lambda,\underline\lambda\}$, we shall study existence of solution, maximum principle and comparison principle, or lack of it, for Dirichlet problems of the following type:
$$
\left\{\begin{array}{lc}
F(x,\nabla u, D^2 u) +b(x)\cdot\nabla u|\nabla u|^{\alpha}
+(c(x)+\lambda)u|u|^\alpha+h(x,u)= f(x)&  \mbox{in } \  \Omega\\
 u=g & \mbox{on }  \partial\Omega.
\end{array}
\right.
$$
depending on the choice of the function $h$.

The maximum principle will be proved under the assumption that 
$h(x,.)$ is non  increasing  and  $h(x,0)=0$.
The comparison principle holds if for all $x\in \Omega$  
$t\mapsto {-h(x,t)\over t^{\alpha+1}}$ is non decreasing on $\R^+$.  We shall also construct a counter-example to 
the comparison principle if this condition fails.

Finally, for continuous functions $f$ and ${\mathcal C}^2$ functions $g$, we shall prove existence of solution
if $h(x,u)=h_1(x,u)- h_2(x,u)$ with 
 $h_i(.,t)\in L^\infty$ for all $t$,
$h_i(x,.)$  non-increasing and continuous for all $x\in\Omega$, $h_i(x,0)=0$ and if
$$\lim_{t\rightarrow\infty} \frac{h_2(x,t)}{t^{\alpha+1}}=0.$$

\section{Notations}

In this section, we state the assumptions on the operators
  
$$G(x,\nabla u,D^2 u): = F(x,\nabla u, D^2 u)+ b(x).\nabla u |\nabla
u|^\alpha $$
treated in this paper,  and the notion of viscosity solution.

The operator $F$ is continuous on $\Omega\times
(\R^N)^\star\times S$, where $S$ denotes the space of symmetric
matrices  on $\R^N$.

The following hypothesis will be considered

\begin{itemize}

\item[(H1)] 
 $F: \Omega\times \R^N\setminus\{0\}\times S\rightarrow\R$, 
and  $\forall t\in \R^\star$,
$\mu\geq 0$,
 $F(x, tp,\mu X)=|t|^{\alpha}\mu F(x, p,X)$.

\item[(H2)] There exist $A\geq a>0$ such that for all $x\in\Omega$, $p\neq 0$ and for all $(M,N)\in S^2$, $N\geq 0$
$$a |p|^\alpha {\rm tr} N\leq F(x,p,M+N)-F(x,p,M)\leq A |p|^\alpha{\rm tr} N$$

\item[(H3)]
There exists a continuous function $\tilde \omega$, $\tilde \omega(0)=0$
such that for all $x,y,$ $p\neq 0$, $\forall X\in S$
$$|F(x,p,X)-F(y,p,X)|\leq \tilde \omega(|x-y|) |p|^\alpha |X|.$$

\item [(H4)]
 There exists a
continuous function
$  \omega$ with $\omega (0) = 0$, such that if $(X,Y)\in S^2$
and 
$\zeta\in \R$ satisfy
$$-\zeta \left(\begin{array}{cc} I&0\\
0&I
\end{array}
\right)\leq \left(\begin{array}{cc}
X&0\\
0&Y
\end{array}\right)\leq 4\zeta \left( \begin{array}{cc}
I&-I\\
-I&I\end{array}\right)$$
where $I$ is the identity matrix in $\R^N$,
then for all  $(x,y)\in \R^N$, $x\neq y$
$$F(x, \zeta(x-y),p, X)-F(y,  \zeta(x-y),p, -Y)\leq \omega
(\zeta|x-y|^2).$$

\end{itemize}

We shall suppose that  $b:\Omega\mapsto
\R^N$ is a continuous and bounded function satisfying:

\begin{itemize}
\item[(H5)]  -Either
$\alpha<0$ and $b$ is H\"olderian of exponent $1+\alpha$,

- or
$\alpha\geq 0$ and, for all $x$ and $y$, 
$$\langle b(x)-b(y), x-y\rangle \leq 0$$
\end{itemize}
\bigskip

Let us recall what we mean by {\it viscosity solutions}, adapted to
our context.

It is well known that, in dealing with viscosity respectively  sub and super
solutions, one works with  
$$u^\star (x) = \limsup_{y, |y-x|\leq r} u(y)$$

and $$u_\star (x) = \liminf_{y, |y-x|\leq r} u(y).$$
 It is easy to see that 
$u_\star \leq u\leq u^\star$ and $u^\star$ is upper semicontinuous (USC)
$u_\star$ is lower semicontinuous (LSC).  See e.g. \cite{CIL, I}.

\begin{definition}\label{def1}

 Let $\Omega$ be a bounded domain in
$\R^N$,  and let $g$ be a given continuous function in $\Omega\times\R$ then
$v$,   bounded on $\overline{\Omega}$, is called a viscosity super solution
of
${G}(x,\grad u,D^2u)=g(x,u)$ if for all $x_0\in \Omega$ if

-Either there exists an open ball $B(x_0,\delta)$, $\delta>0$  in $\Omega$
on which 
$v= cte= c
$ and 
$0\leq g(x,c)$, for all $x\in B(x_0,\delta)$

-Or
 $\forall \varphi\in {\mathcal C}^2(\Omega)$, such that
$v_\star-\varphi$ has a local minimum on $x_0$ and $\grad\varphi(x_0)\neq
0$, one has
\begin{equation}
{G}( x_0,\grad\varphi(x_0),
 D^2\varphi(x_0))\leq g(x_0,v_\star(x_0)).
\label{eq1}\end{equation}

Of course $u$ is a viscosity sub solution if
for all $x_0\in \Omega$,

-Either there exists a ball $B(x_0, \delta)$, $\delta>0$ on which  $u =
cte= c$ and
$0\geq g(x,c)$, for all $x\in B(x_0,\delta)$

-Or  $\forall
\varphi\in {\mathcal C}^2(\Omega)$, such that
$u^\star-\varphi$ has a local maximum on $x_0$ and
$\grad\varphi(x_0)\neq 0$, one has
\begin{equation}
 G( x_0,\nabla \varphi(x_0),
D^2\varphi(x_0))\geq g(x_0,u^\star(x_0)) . \label{eq2}\end{equation}

A	 viscosity solution is a function which is both a super-solution and a sub-solution.
\end{definition}

See e.g. \cite{CGG} for a similar definition of viscosity solution of 
equations with singular operators.

For convenience we recall  the definition of semi-jets given e.g. in 
\cite{CIL}

\begin{eqnarray*}
J^{2,+}u(\bar x) &= &\{ (p, X)\in \R^N\times S, \ u(x)\leq
u(\bar x)+
\langle p, x-\bar x\rangle + \\
& &+{1\over 2} \langle X(x-\bar x), x-\bar
x\rangle
+ o(|x-\bar x|^2) \}
\end{eqnarray*}
and 

\begin{eqnarray*}J^{2,-}u(\bar x)& =& \{ (p, X)\in \R^N\times S, \
u(x)\geq u(\bar x)+
\langle p, x-\bar x\rangle +\\
&&+ {1\over 2} \langle X(x-\bar x), x-\bar
x\rangle
+ o(|x-\bar x|^2\}.
\end{eqnarray*}
In the definition of viscosity solutions the test functions can be
substituted by the elements of the semi-jets in the sense that in the
definition above one can restrict to the functions $\phi$ defined by
$\phi(x)=u(\bar x)+
\langle p, x-\bar x\rangle 
+ {1\over 2} \langle X(x-\bar x), x-\bar
x\rangle$  with $(p,X)\in J^{2,-}u(\bar x)$ when $u$ is a super solution, 
and $(p,X) \in J^{2,+}u(\bar x)$ when $u$ is a sub solution.  

\bigskip

For the convenience of the reader we now recall the properties obtained for $\lambda<\overline\lambda$ in \cite{BD2}.

\begin{theorem}\label{maxp}

Let $\Omega$ be a bounded domain of $\R^n$.
Suppose that $F$ satisfies (H1), (H2), (H4),  that $b$ and $c$ are continuous and
$b$ satisfies (H5).  Suppose that
$\tau<\bar\lambda$ and that
$u$ is a  viscosity sub solution of 

$$G(x,\nabla u, D^2 u)+ (c(x)+ \tau) |u|^{\alpha}u \geq 0\  \ {\rm in}\ \ \Omega$$
with $u\leq 0$ on the boundary of $\Omega$, then 
$u\leq 0$ in $\Omega$. 

If $\tau<
\underline\lambda$ and  $v$ is a super solution of 
$$G(x, \nabla v, D^2 v)+(c(x)+ \tau) |v|^{\alpha}v \leq 0\  \ {\rm in}\ \ \Omega$$
with $v\geq 0$ on the boundary of $\Omega$ then $v\geq 0$ in $\Omega$.

\end{theorem}

Let us also recall the following comparison theorem.
\begin{theorem}\label{complambda1}
Suppose that $F$ satisfies (H1), (H2), and (H4),  
that $b$ and $c$ are continuous and
bounded and $b$ satisfies
$(H5)$.  Suppose that $\tau< \bar\lambda$, $f\leq 0$, $f$ is
upper semi-continuous and
$g$ is lower semi-continuous  with $f\leq  g$.

Suppose that there exist
$\sigma$ upper semi continuous , and
$v$ non-negative and lower semi continuous , satisfying
\begin{eqnarray*}
G(x, \grad v,D^2v)+(c(x)+\tau) v^{1+\alpha} & \leq & f \quad 
\mbox{in}\quad \Omega \\ G(x,  \grad \sigma,D^2\sigma)
+(c(x)+\tau)|\sigma|^{\alpha}\sigma & \geq & g \quad  \mbox{in}\quad
\Omega  \\ 
\sigma \leq  v &&   \quad  \mbox{on}\quad \partial\Omega 
\end{eqnarray*}
Then $\sigma\leq v$ in $\Omega$ in each of these two cases:

\noindent 1) If $v>0$ on $\overline{ \Omega}$ and either $f<0$ in 
$\Omega$, 
 or  $g(\bar x)>0$ on every point $\bar x$ such that  $f(\bar x)=0$,

\noindent  2) If $v>0$ in $\Omega$, $f<0$ and $f<g $ on
$\overline\Omega$ 
\end{theorem}

We recall some of the properties of the distance function for bounded
${\mathcal C}^2$ set. 
\begin{remark}\label{dist} In all the paper we shall consider that $\Omega$ is a bounded $\mathcal C^2$ domain. In particular we shall use several times the fact that this implies  that the distance to the boundary:
$$d(x,\partial\Omega):=d(x):=\inf\{|x-y|,\  y\in\partial\Omega\}$$
satisfies the following properties:
\begin{enumerate}
\item $d$ is Lipschitz continuous
\item There exists $\delta>0$ such that in $\Omega_\delta=\{x\in\Omega\ \mbox{such that }\ d(x)\leq \delta\}$, $d$ is $\mathcal C^{1,1}$.
\item $d$ is semi-concave, i.e. there exists $C_1>0$ such that 
$d(x)-C_1|x|^2$ is concave and this implies $J^{2,+} d(x)\neq \emptyset$ for any $x\in\Omega$.
\item If  $J^{2,-}d(x)\neq \emptyset$, $d$ is  differentiable at $x$
and
$|\nabla d(x)|=1$. 
\end{enumerate}
\end{remark}

\bigskip

\section{An existence's result for ${\mathcal C}^2$ boundary data  }

 We assume that $F$ satisfies the assumptions enumerated in the previous
section.

Let $g$ be given in $W^{2,\infty}(\partial \Omega)$, and $f\in L^\infty$.
We denote by $b$ a bounded and continuous function which satisfies 
(H5).
We introduce some  functions $h_1$ and $h_2$ :  

\bigskip
(H6) For $i=1$ and $i=2$ let $h_i:\Omega\times \R\rightarrow \R$ such that $h_i(.,t)\in L^\infty$ for all $t$,
$h_i(x,.)$ is non-increasing and continuous for all $x\in\Omega$, $h_i(x,0)=0$ and
$$\lim_{t\rightarrow\infty} \frac{h_2(x,t)}{t^{\alpha+1}}=0.$$

\bigskip
Our main existence's result is the following theorem:

\begin{theorem}\label{mexi2}
Suppose that $\lambda < \lambda_1=\inf\{\overline{\lambda},
\underline{\lambda}\} $ and suppose that $h_1(x,t)$ and $h_2(x,t)$
satisfy (H6) then, for $g\in W^{2,\infty} (\partial \Omega)$ and $f\in
{\mathcal C}(\overline{\Omega})$, there exists a solution of 
$$
\left\{\begin{array}{lc}
G(x,\nabla u, D^2 u) 
+(c(x)+\lambda)u|u|^\alpha+h_1(x,u)= f(x)+h_2(x,u) &  \mbox{in } \  \Omega\\ u=g & \mbox{on }
\partial\Omega.
\end{array}
\right.
$$
\end{theorem}

In particular, Theorem \ref{mexi2} implies that, for 
$\beta_1(x)\geq 0$ and $\beta_2(x)\geq 0$,  for any $q_1>0$ and for 
$q_2<\alpha$, there exists $u$ solution of
$$
\left\{\begin{array}{rc}
G(x,\nabla u, D^2 u) 
+(c(x)+\lambda)|u|^\alpha u-\beta_1(x)|u|^{q_1}u = & \\
 f(x)-\beta_2(x)|u|^{q_2}u &  \mbox{in } \  \Omega\\
 u=g & \mbox{on }
\partial\Omega.
\end{array}
\right.
$$

\bigskip
\begin{remark} \label{remexis} The previous existence's result still holds when $\underline\lambda<\lambda<\overline\lambda$ if $f\leq 0$ and $g\geq 0$.
The proof proceeds as the one of Theorem \ref{mexi2} using the Remark \ref{remconst}
which is stated after Theorem \ref{const} and the fact that $u\equiv 0$ is a subsolution.

A symmetric result holds for $\underline\lambda>\lambda>\overline\lambda$.
\end{remark}

The proof of this Theorem requires  several step : 
 
 The first step is given by:
 \begin{proposition}\label{prop0}

Suppose that $g$ is in  $W^{2,\infty} (\partial \Omega)$, that
$h:\Omega\times\R\rightarrow \R$ is such that $h(x,.)$ is non increasing
and continuous and $f$ is in $L^\infty$. Then there exists $u$ a
viscosity solution of 
\begin{equation}\label{eqprop0}
\left\{
\begin{array}{lc}
G(x,\nabla u, D^2 u)+h(x,u)= f &  
\mbox{in}\ \Omega   \\ u=g  &    \mbox{on}\ \partial\Omega.   \\
 \end{array}
\right.
\end{equation}
\end{proposition}
 
 To prove this proposition,  it is enough to construct a sub and a super solution of
 (\ref{eqprop0})
and then  apply Perron's method (see  \cite{BD7, I}).
This is the purpose of the following two Propositions \ref{prop1}, \ref{prop2}.

Then, in Theorem \ref{Hol},  we will prove a H\"older's estimate  which also 
gives a compactness result. And then the proof of  Theorem \ref{mexi2}
will be done through a recursive argument.

\begin{proposition}\label{prop1}

Suppose that $g$ is  in $W^{2,\infty}(\partial \Omega)$, that
$h:\Omega\times\R\rightarrow
\R$ is such that $h(x,.)$ is non increasing and continuous and $m\in
\R^+$. Then there exists $u$ a viscosity subsolution of 
\[
\left\{
\begin{array}{lc}
G(x,\nabla u, D^2 u)+h(x,u)\geq m &  
\mbox{in}\ \Omega   \\ u=g  &    \mbox{on}\ \partial\Omega.   \\
 \end{array}
\right.
\]
\end{proposition}

{\bf Proof of Proposition \ref{prop1}.}
We denote by $g$  a  ${\mathcal C}^2\cap W^{2,\infty} (\Omega)$
function which equals 
$g$ on $\partial\Omega$.

Let us note first that it is sufficient to construct a subsolution $v$ when $h\equiv 0$, as long as it satisfies
$v\leq g$.

 Indeed, suppose that we have constructed  a subsolution  of the equation

$$G(x,\nabla v, D^2 v) \geq m^\prime:= m-h(x,g) $$
this will imply, since $h$ is decreasing, that 
$$G(x,\nabla v, D^2 v)+h(x,v)\geq
m-h(x,g)+h(x,v)\geq m.$$

We  assume now that $h\equiv 0$.

 Let us consider, for $L$ and $k$ to be chosen later, the function
$v(x) = g(x)+L ((1+d(x))^{-k}-1)$.
Let $D$ be an upper bound for $d$ on $\Omega$, and  
$C_1$ be a  positive constant  as in 3 of Remark \ref{dist}.

We choose $k>1 $ large enough in order that 
$a(k+1)> 2(A+a) NC_1(1+D)$ and 
$k+1> \frac{2^{\alpha+2}}{a}|b|_\infty (1+D)^{\alpha+2}$. Next we choose  $L$ such that 
$$L k\geq  \sup \left\{2|\nabla g |_\infty (1+D)^{k+1}, 2{(A+a)\over a}
(1+D)^{k+2}|D^2 g|_\infty, ({16m
^\prime (1+D)^{k(\alpha+1)+\alpha+2}
\over a (k+1)})^{1\over 1+\alpha}\right\}.$$

Let us recall that since $d$ is semi-concave for every $x\in \Omega$ 
$J^{2,+} d(x)\neq \emptyset$. 

To prove that $v$ is a subsolution, let $x_0\in \Omega $
 be any point such that
there exists a test function $\varphi\in C^2$  satisfying
$$g(x)+({1\over (1+d)^k}-1)(x)\leq \varphi(x),\quad g(x_0)+ ({1\over (1+d)^k}-1)(x_0)=
\varphi(x_0). $$
 Then ${1\over (\varphi(x)-g(x)+1)^{1\over k}}-1\leq d(x)$ 
and 
$J^{2,-}d(x_0)\neq \emptyset$ .

Since $J^{2,+} d(x_0)\neq \emptyset$  this implies that $d$ is
differentiable in
$x_0$ and then  $\grad d(x_0)=\grad \left({1\over (\varphi-g+1)^{1\over k}}-1\right)(x_0)$ which of course implies that 
$$\grad v(x_0)=\grad\varphi(x_0)= (\nabla g -L k(1+d)^{-k-1}\grad d)(x_0).$$

Hence, for simplicity, we shall use directly  $v$ instead of the test 
function. Moreover  since $d$ is semiconcave, one has, for the constant $C_1$ defined in Remark \ref{dist}, for all $x$ and
$\bar x$ 
$$d(x)-C_1|x|^2 -d(\bar x)+ C_1 |\bar x|^2 \leq \nabla d(\bar x)-C_1 \bar
x. (x-\bar x).$$
This implies that for every $(p,X)\in J^{2,- } d(\bar x)$,
$X\leq C_1I $.

 Observe that in
particular, with the above  choice of
$L$ and
$k$ one has 
$$|\nabla v|\geq {L k(1+D)^{-k-1}\over 2}.$$

We are now in a position to compute $D^2v$:

$$D^2 v = D^2 g + L
k(k+1)(1+d)^{-k-2}\nabla d\otimes \nabla d-L k(1+d)^{-k-1} D^2 d.$$
Recalling that $\grad d\otimes\grad d\geq 0$ and $|\grad d|=1$, we get

\begin{eqnarray*}
{\mathcal M}_{a,A}^-(D^2v) &\geq& Lk \left(a{(k+1)\over
(1+d)^{k+2}}-{(A+a) C_1N\over (1+d)^{k+1}}\right)-(A+a) |D^2
g|_\infty\\
&\geq& {Lk a(k+1)\over 2(1+d)^{k+2}}-(A+a) |D^2 g|_\infty\\
&\geq &{Lk(k+1)a\over 4 (1+d)^{k+2}}.
\end{eqnarray*}
As a consequence 
\begin{eqnarray*}
|\nabla v|^\alpha {\mathcal M}_{a,A}^-(D^2v) + b(x).\nabla
v|\nabla v|^\alpha
&\geq &{(Lk)^{\alpha+1}\over (1+D)^{k(\alpha+1)+\alpha+2}} \left(
{a(k+1)\over 2^{\alpha+2}}-|b|_\infty (1+D)^{\alpha+2}\right)\\
&\geq& 
{(Lk)^{\alpha+1}a(k+1)\over 2^{\alpha+3}(1+d)^{k(\alpha+1)+\alpha+2}}\\
&\geq&
2m^\prime .
\end{eqnarray*}

We then obtain the required inequality,
\begin{eqnarray*}
G(x,\nabla v, D^2 v) &\geq&
{\mathcal M}_{a,A}^-(D^2v) |\nabla v|^\alpha -|b|_\infty |\nabla v|^{\alpha+1} \\
&\geq& m^\prime
\end{eqnarray*}
This ends the proof.

\begin{proposition}\label{prop2}
Suppose that $m>0$ that  $g\in W^{2,\infty} (\partialÊ\Omega)$ and
$h:\Omega\times\R\rightarrow \R$ is such that $h(x,.)$ is non increasing
and continuous . Then there exists a supersolution $u$ of
$$\left\{\begin{array}{lc}
F(x,\nabla u, D^2 u)+b(x).\nabla u |\nabla u|^\alpha+h(x,u)  \leq -m
&
\mbox{in }\
\Omega\\ u = g&   \mbox{ on }  \partial \Omega.
\end{array}
\right.
$$ 
\end{proposition}
\begin{remark}
In the following we shall denote by $\overline{u}=S(g,-m) $ a
supersolution as in    Proposition \ref{prop2} and by
$\underline{u} = S(g,m)$ a subsolution as in Proposition \ref{prop1}.
\end{remark}
{\bf Proof of Proposition \ref{prop2}.}
We still denote by $g$ a ${\mathcal C}^2\cap W^{2,\infty}$ function on
$\overline{\Omega}$ which equals $g$ on the boundary. 
As in the previous proposition it is sufficient to prove the result when
$h\equiv0$ as long as
$\varphi\geq g$.

We choose  
  $\varphi(x):= g(x)+ L
(1-(1+d(x))^{-k})$, with $L$ and $K$ appropriate constants to be chosen later.

To prove that $\varphi$ is a supersolution, either $J^{2,-}
\varphi(x_0)= \emptyset$ or not, and then there exists $\psi(x)\leq g(x)+
L(1-(1+d)^{-k})(x)$ touching $\varphi$ at $x_0$.
This implies that 

$1-(1-g(x)+\varphi)(x)^{-k}(x)\geq d(x)$, and  $J^{2,-}d(x)\neq \emptyset.$
Then $J^{2,+}d(x_0)$ and $J^{2,-}d(x_0)$ are non empty, and $d$ is
differentiable on $x_0$. We shall use in the following $\nabla d$  for
the computations below. 
 As in the previous proof, $k$ is chosen large
enough in order that

$$(k+1)\geq \sup \left(\frac{2(A+a)C_1N(1+D)}{a}, 
\frac{|b|_\infty4(1+D)}{a}\right).$$ Now  we can choose $L$ such that
$$L\geq
\sup\left(\left(\frac{4m'}{a(k+1)}\right)^{\frac{1}{\alpha+1}}
\frac{2(1+d)^{k+1}}{k} 
, \frac{2(A+a)|D^2g|_\infty(1+D)^{k+2}}{ka(k+1)},
\frac{2|\grad g|_\infty}{k(1+D)^{k+1}}\right).$$

The computation of the gradient  gives 
$$\nabla \varphi = \nabla g+L k 
(1+d)^{-k-1} \nabla d.$$
and by the previous assumptions $|\nabla \varphi|\geq {
L\over 2 }k 
(1+d)^{-k-1} $.

While
$$D^2\varphi= D^2g+\frac{Lk}{(1+d)^{k+1}}\left[-\frac{(k+1)}{(1+d)} \grad d\otimes\grad d +d D^2d\right].$$

We then have that

$${\mathcal M}_{a,A}^+(\nabla \nabla \varphi) \leq (A+a)|D^2g|_\infty
+\frac{Lk}{(1+d)^{k+1}}\left(-a\frac{(k+1)}{(1+d)} +(A+a)dC_1N\right).$$
Since
$$
a\frac{(k+1)}{2(1+d)}\geq (A+a)d C_1N,$$
one gets,

$${\mathcal M}_{a,A}^+(\nabla \nabla \varphi) \leq (A+a)|D^2g| -\frac{Lk(k+1)a}{2(1+d)^{k+2}} \leq
-\frac{Lka(k+1)}{4(1+d)^{k+2}}.$$

We can finish the computation, and get;
\begin{eqnarray*}
|\nabla \varphi|^\alpha {\mathcal M}_{a,A} (\nabla \nabla
\varphi) + b(x).\nabla \varphi|\nabla \varphi|^\alpha&\leq&
\left({Lk\over 2(1+d)^{k+1}}\right)^{\alpha +1}\left(-\frac{a(k+1)}{2(1+d)}+|b|_\infty\right)
\end{eqnarray*}

Hence, with our choice of the constants, we have obtained   
$$
G(x,\nabla \varphi, D^2\varphi)  
\leq  -m^\prime.
$$
This ends the proof.
 
\bigskip
We now give regularity results for solutions with 
boundary data in  $W^{2,\infty}$.

This immediately implies a compacity result for fixed data
$W^{2,\infty}(\partial \Omega)$,  extending in that way the results in
\cite{BD7}.
\begin{theorem}\label{Hol}

Suppose that $g\in W^{2,\infty}(\partial\Omega)$  
then every solution of 
$$
\left\{\begin{array}{lc}
G(x,\nabla u, D^2 u)= f & \mbox{in}\
\Omega\\ u=g& \mbox{on }\ \partial\Omega
\end{array}
\right.
$$
satisfies: For any $\gamma\in(0,1)$ there exists $C$ depending on 
$\gamma$, $|g|_{ W^{2,\infty} (\partial \Omega)}$ and $|f|_\infty$
$$ |u(x)-u(y)|\leq C|x-y|^{\gamma}$$
for all $x$ and $y$ in $\Omega$.
\end{theorem}
\begin{corollary}\label{comp}
Suppose that $g_n\rightarrow
 g$ in $W^{2,\infty} (\partialÊ\Omega)$, that $f_n$ is a bounded sequence  of bounded functions and
that $u_n$ are  solutions of 

$$\left\{
\begin{array}{lc}
G(x,\nabla u_n, D^2u_n) = f_n & \ {\rm in } \ \Omega\\
u_n=g_n & \ {\rm on} \ \partial\Omega
\end{array}
\right.
$$
Suppose that $(u_n)$ is bounded in $L^\infty$, then $u_n$ is uniformly
H\"olderian, and  the sequence  is relatively compact in
${\mathcal C} (\overline{\Omega})$.
\end{corollary}
{\bf Proof of Theorem \ref{Hol} .} We proceed similarly to \cite{BD7,IL}. Let us recall that
the proof has two steps. In the first one H\"older's regularity is proved near the boundary.
And then it is proved in the interior, through a
 typical viscosity argument.
 We only give the details of the 
first part, since the second part proceeds as in \cite{BD7}.

Hence we define for a given positive $\delta$ 
$$\Omega_\delta:=\{(x,y)\in\Omega^2\  \mbox{such that} \  |x-y|\leq\delta\}$$
and the first step consists in proving that
on $\partial\Omega_\delta$ there exists $C>0$ such that
$$u(x)-u(y)\leq C|x-y|^\gamma.$$

If $C> \frac{2|u|_\infty}{\delta}$ then the inequality is true for $|x-y|=\delta$ so we should only prove it for
$(x,y)\in\Omega\times\partial\Omega$ and similarly for $(x,y)$  $\in
\partial \Omega\times \Omega$.

We shall still denote by $g$ a $C^2\cap W^{2,\infty } (\Omega)$ extension
of
$g$ to
$\Omega$.

Using the Propositions \ref{prop1} and \ref{prop2} we know that there exist $L$, $k$,  which depend
only on  universal constants and on the $C^2$ norm of $g$, such
that 
$g(x)+L ((1+d)^{-k}-1)$ is  a subsolution and $g(x)+ L
(1-(1+d)^{-k})$ is a super solution.

Using the comparison theorem in \cite{BD1}, one gets that  
$$g(x)+L ((1+d)^{-k}-1)\leq u \leq g(x)+ L
(1-(1+d)^{-k})$$ 

Finally there exist $C$ and $C_g$ some Lipschitz constant of $g$ such
that  on $\partial \Omega_\delta $,
if $d(x) = \delta$ and $y\in \partial \Omega$, one has 
$$u(x)-u(y) = u(x)-g(y) \geq  g(x)-g(y)-Cd(x)^\gamma \geq
-C_g |x-y|-C|x-y|^\gamma$$
and 
$$u(x)-u(y) \leq u(x)-g(y)\leq g(x)-g(y)+C(d(x))^\gamma
\leq C_g|x-y|+ C |x-y|^\gamma.$$
The rest of the proof proceeds as in \cite {BD7}.

\bigskip
We pass to the second step, which treats the case where $f$ and $g$ are
constant, with opposite sign.  
\begin{theorem}\label{const}
Suppose that we have the same hypothesis as in Theorem \ref{mexi2} in particular for $\lambda<\lambda_1$.
Then for all $g\in \R^+$, $f\in \R^+$ there exist $\overline{u}\geq 0$ and
$\underline{u}\leq 0$ respectively solutions of 
$$
\left\{\begin{array}{lc}
G(x,\nabla \overline{u} , D^2\overline{u})+ (c(x)+\lambda)
|\overline{u}|^\alpha \overline{u}+h_1(x,\overline u)= -f+h_2(x,\overline
u) &
\mbox{in}\
\Omega\\ \overline{u}=g& \mbox{on }\ \partial\Omega
\end{array}
\right.
$$
and 
$$
\left\{\begin{array}{lc}G(x,\nabla \underline{u}, D^2\underline{u})+
(c(x)+\lambda) |\underline{u}|^\alpha \underline{u} +h_1(x,\underline u)=
f + h_2(x,\underline u)&
\mbox{in}\
\Omega\\ \underline u=-g& \mbox{on }\ \partial\Omega.
\end{array}
\right.
$$
\end{theorem}
\begin{remark}\label{remconst}
If $\underline\lambda<\lambda<\overline\lambda$ the same existence of $\overline u$ holds when $f\leq 0$ and $g\geq 0$
and symmetrically if  $\underline\lambda>\lambda>\overline\lambda$ for $\underline u$.
\end{remark}

{\bf Proof of Theorem \ref{const}.}
We consider the first case, the other being symmetric.

Let $u_n$ be the  sequence defined as the positive solution of 
$$\left\{\begin{array}{rc}
G(x, \nabla u_n, D^2 u_n)
(c+\lambda-2|c+\lambda|_\infty)u_n^{1+\alpha}+ h_1(x,u_n) = & \\
-f-2|\lambda+c|_\infty
u_{n-1}^{1+\alpha} +h_2(x,u_{n-1})& \mbox{in}\  \Omega\\ u_n = g, & \mbox{on
}\partial\Omega\\
\end{array}
\right.
$$
for $n\geq 1$ and $u_0=0$, which exists by Proposition \ref{prop0}.

By construction and the comparison principle, $u_n$ is 
increasing and hence $u_n\geq 0$. Let us prove that it is bounded. Indeed if not,
$|u_n|_\infty\rightarrow+\infty $ and  the function $w_n= {u_n\over
|u_n|_\infty}$ solves
\begin{eqnarray*}
&&G(x, \nabla w_n, D^2 w_n)+
 (c+\lambda-2|c+\lambda|_\infty)w_n^{1+\alpha}+
\frac{h_1(x,u_n)}{|u_n|_\infty^{\alpha+1}}   \\
 &=& {f\over |u_n|_\infty^{1+\alpha}}-
2|\lambda+c|_\infty \frac{u_{n-1}^{1+\alpha}}{|u_n|_\infty^{1+\alpha}}+
\frac{h_2(x,u_{n-1})}{|u_n|_\infty^{\alpha+1}}.
\end{eqnarray*}
By hypothesis on $h_1$ and $h_2$, 
$h_1(x,u_n)\leq 0$ and
$$\frac{h_2(x,u_{n-1})}{|u_n|_\infty^{\alpha+1}}\geq \frac{h_2(x,|u_n|_\infty)}{ |u_n|_\infty^{\alpha+1}}\rightarrow 0.$$

 Extracting subsequences one has for some increasing sequence $\sigma(n)$ 
$${|u_{\sigma(n)-1}|_\infty\over |u_{\sigma(n)}|_\infty}\rightarrow k\leq 1$$
and 
$$w_{\sigma(n)}\rightarrow w\geq 0$$
with $|w|_\infty=1$ and $w=0$ on the boundary of $\Omega$. Moreover, using the
compactness result in Corollary \ref{comp}, one gets that $w$ satisfies 
$$G(x,\nabla w, D^2 w)+(c+\lambda-2|c+\lambda|_\infty(1-k^{\alpha+1}))
w^{1+\alpha}\geq 0$$ with $w=0$ on the boundary. Since
$\lambda-2|c+\lambda|_\infty(1-k^{\alpha+1})< \lambda_1$  , by the maximum
principle we get that $w\leq 0$. But $w\geq 0$  and hence $w\equiv 0$ ,
 which contradicts $|w|_\infty = 1$. 

We have proved that $(u_n)$ is bounded. Since $(u_n)$ is monotone, by the compactness result of
Corollary
\ref{comp},  it
converges to some $u$ which is the required solution.
This ends the proof.

\bigskip
{\bf Proof of Theorem \ref{mexi2}.}
We need to construct a solution of 
$$\left\{\begin{array}{lc}
G(x,\nabla u, D^2 u)+
(c+\lambda)|u|^{\alpha}u+ h_1(x,u) = f+h_2(x,u) & \mbox{in}\  \Omega\\ u = g, & \mbox{on
}\ \partial\Omega\\
\end{array}
\right.
$$
with $f$ in $L^\infty$ and $g$ which is ${\mathcal C}^2$. 

We denote respectively  $\underline{u}=S(2|f|_\infty,-|g|_\infty)$  $\overline{u}=S(-2|f|_\infty,|g|_\infty)$ the solutions obtained in the
previous theorem.

We define a sequence $u_n$ by the recursive process  :
$$\left\{\begin{array}{rc}
G(x,\nabla u_n, D^2 u_n)+
((c+\lambda)-2|c+\lambda|_\infty )|u_n|^{\alpha}u_n +h_1(x,u_n)& \\
 =
f-2|c+\lambda|_\infty |u_{n-1}|^\alpha u_{n-1} +h_2(x,u_{n-1}) &
\mbox{in}\
\Omega\\ u_n = g, &
\mbox{on }\partial\Omega\\
\end{array}
\right.
$$
 initializing with  $u_0 = \underline{u}$.  We know that the sequence is well defined by Proposition \ref{prop0} with $h(x,u)=((c+\lambda)-|c+\lambda|_\infty )|u|^{\alpha}u+h_1(x,u)$.

Let us prove that 
$$\underline{u}\leq u_n\leq \overline{u}.$$
 Let us note that since
$\underline{u}\leq 0\leq \overline{u}$ it is equivalent to prove that 
$u_n^+\leq \overline{u}$ and $u_n^-\leq -\underline{u}$.

Suppose that we know that $\underline{u}\leq u_{n-1}\leq \overline{u}$
and let us prove that 
$\underline{u}\leq u_{n}\leq \overline{u}$.

One has 
\begin{eqnarray*}
&& G(x,\nabla u_n, D^2u_n) 
+(c+\lambda-2|c+\lambda|_\infty)u_n|u_n|^{\alpha}+h_1(x,u_n)\\
&= & f-2|\lambda+c|_\infty
|u_{n-1}|^\alpha u_{n-1}+h_2(x,u_{n-1})\\
&\leq &|f|_\infty +2|\lambda+c|_\infty
|u_{n-1}^-|^{1+\alpha}+h_2(x,-u_{n-1}^-) \\
&\leq &|f|_\infty -2| \lambda+c|_\infty
|\underline{u}|^\alpha\underline{u}+h_2(x,\underline u)\\
&< &G(x,\nabla \underline{u}, D^2 \underline{u} )
+(c+\lambda-2|c+\lambda|_\infty)\underline u|\underline
u|^{\alpha}+h_1(x,\underline u).
\end{eqnarray*}
Since on the boundary $u_n\geq \underline{u}$  one gets that 
$u_n\geq \underline{u}$.

Now 
\begin{eqnarray*}
&& G(x,\nabla u_n, D^2 u_n) 
+(c+\lambda-2|c+\lambda|_\infty)u_n|u_n|^{\alpha} +h_1(x,u_n)=\\
& = & f-2|\lambda+c|_\infty
|u_{n-1}|^\alpha u_{n-1}+h_2(x,u_{n-1})\\
&\geq &-|f|_\infty -2|\lambda+c|_\infty
|u_{n-1}^+|^{1+\alpha} +h_2(x,u_{n-1}^+)\\
&\geq& -|f|_\infty -2| \lambda+c|_\infty
\overline{u}^{\alpha+1}+h_2(x,\overline u) \\
& >& G(x,\grad\overline{u}, D^2\overline{u})
+(c+\lambda-2|c+\lambda|_\infty)\overline
u^{\alpha+1}+h_1(x,\overline u).
\end{eqnarray*}
 Since on the boundary 
$u_n= g\leq \overline{u}$ one gets that 
$$u_n\leq \overline{u}.$$

It is sufficient now to invoke the compacity result to see that the
sequence $u_n$
 is relatively compact in ${\mathcal C} (\overline{\Omega})$. One gets, 
passing to the limit,  that for a subsequence $u_n\rightarrow u$ which satisfies
$$\left\{\begin{array}{lc}
G(x,\nabla u, D^2 u)+
(c(x)+\lambda)|u|^{\alpha}u +h_1(x,u) = f+h_2(x,u) & \mbox{in} \  \Omega\\ u = g, & \mbox{on
}\ \partial\Omega.\\
\end{array}
\right.
$$  
This ends the proof.

\subsection{The case $f(x)\equiv0$.}

\begin{proposition} Suppose that for all $x$ and $u$, $h$ satisfies
$h(x,u)u\leq 0$. Then for
$\lambda<\lambda_1$ the only solution of
$$
\left\{\begin{array}{lc}
G(x,\nabla u, D^2 u)+ (\lambda+c)
|u|^\alpha u +h(x,u)=0 & \mbox{in }\ \Omega \\ u=0,& \mbox{on
}\partial\Omega\\
\end{array}
\right.
$$
is $u\equiv 0$.

\end{proposition}
{\bf Proof}
Suppose that there exists a non zero solution.
Then either $\Omega^+ = \{ x, u(x)>0\}\neq\emptyset$ or $\Omega^-
 = \{ x, u(x)<0\}\neq \emptyset$. Without loss of generality we can
suppose that  $\Omega^+\neq\emptyset$
 then
$$G(x,\nabla u, D^2u)+ (\lambda+c)
|u|^\alpha u = -h(x,u)\geq 0$$ in $\Omega^+$
and $u=0$ on the boundary of $\Omega^+$. 
As seen in \cite{BD7}, $\overline\lambda(\Omega^+)
\geq \overline\lambda(\Omega)>\lambda$, then by the maximum principle 
$u\leq 0$ in $\Omega^+$ which is a contradiction. 
This ends the proof.

\bigskip

We consider here $\beta$ some continuous and bounded function and
$q<\alpha$. 

We are interested in the existence of non trivial solutions of

\begin{equation}
\left\{\begin{array}{lc}\label{eqnonh}
G(x,\nabla u, D^2 u)+ (\lambda+c)
|u|^\alpha u +\beta |u|^q u=0 & \mbox{in }\ \Omega \\ u=0,& \mbox{on
}\partial\Omega\\
\end{array}
\right.
\end{equation}

\begin{theorem}\label{nonh}
Suppose that $\lambda<\lambda_1$. Then, if $\beta=\beta^+-\beta^-$ with $\beta^+\not\equiv 0$ there exists a non trivial solution of (\ref{eqnonh}).
\end{theorem}

\bigskip

\noindent
{\bf Proof of Theorem \ref{nonh}.}
First we suppose that $\beta\geq 0$, $\beta$  not identically zero.
We begin to construct a subsolution with the aid of some eigenfunction.
Let $\phi>0$ be such that 
$$
\left\{ \begin{array}{lc}G(x, \nabla \phi, D^2\phi)+(\overline{\lambda} +c)
\phi^{1+\alpha} = 0& \mbox{in}\ \Omega\\
\phi=0 & \mbox{on}\  \partial \Omega,
\end{array}
\right.$$
with $|\phi|_\infty=1$.
If $m$ is small enough in order that 
$$(m\phi)^{\alpha-q}|\beta|_\infty < {\overline{\lambda}-\lambda\over 2}$$
then $m\phi$ is a subsolution of (\ref{eqnonh}). 

Let $u_n$ be defined in a  recursive way by 
$$\left\{\begin{array}{cc}
G(x, \nabla u_n, D^2 u_n)+(c+\lambda-2|c+\lambda|_\infty|)
|u_n|^{\alpha} u_n = & \\ = -2|\lambda+c|_\infty |u_{n-1}|^\alpha
u_{n-1}-\beta |u_{n-1}|^q u_{n-1},&\mbox{in}\ \Omega
\\ u_n=0,& \mbox{on
}\partial\Omega.\\
\end{array}
\right.
$$ 
with $u_0=m\phi$.
The solutions $u_n$ are well defined by Proposition \ref{prop0}. We begin by proving that $u_1\geq m\phi$.

\begin{eqnarray*}
&& G(x,\nabla u_1, D^2 u_1)+(c+\lambda-2|c+\lambda|_\infty|)
|u_1|^{\alpha} u_1\\ &=& -2|\lambda+c|_\infty u_{0}^{\alpha+1}-\beta
u_{0}^{q+1}\\ &\leq &-2|\lambda+c|_\infty
m^{1+\alpha}\phi^{\alpha+1}-m^{q+1}\beta
\phi^{q+1}\\
&\leq & G(x, \nabla (m\phi), D^2(m\phi)) +(c+\lambda-2|c+\lambda|_\infty|)
(m\phi)^{\alpha+1}
\end{eqnarray*}
this implies, by the comparison principle, that 
$u_1\geq m\phi$. 

The same reasoning establishes that $u_n\geq u_{n-1}$. 
Suppose that $|u_n|_\infty\rightarrow +\infty$,
then, defining $w_n = {u_n\over |u_n|_\infty}$, one gets that 
$$\left\{\begin{array}{rc}
G(x,\nabla w_n D^2 w_n) +(c+\lambda-2|c+\lambda|_\infty) |w_n|^{\alpha}
w_n= & \\
-2|\lambda+c|_\infty
|w_{n-1}|^\alpha w_{n-1}{|u_{n-1}|^{1+\alpha}_\infty \over
|u_n|_\infty^{1+\alpha}}-\beta {|u_{n-1}|_\infty^q\over |u_n|_\infty^{1+\alpha}}|w_{n-1}|^q
w_{n-1},&\mbox{in}\
\Omega
\\ w_n=0,& \mbox{on
}\partial\Omega.\\
\end{array}
\right.
$$  

By the compactness result,  $w_n$  converges, up to a subsequence, to  some non negative function
$w$ which is such that $|w|_\infty = 1$, and for some $k\leq 1$ it is a solution of
$$\left\{\begin{array}{lc}
G(x,\nabla w, D^2w) + (c+\lambda-2|c+\lambda|_\infty(1-k^{\alpha+1}))|w|^{1+\alpha} w=0 & \mbox{in}\ \Omega \\
w=0 & \mbox{on}\ \partial \Omega.
\end{array}
\right.$$
Then, since $\lambda -2|\lambda+c|(1-k^{1+\alpha})<
\overline{\lambda}$, one gets that $w=0$, a contradiction. 
Finally the sequence $u_n$ is increasing and bounded and by the
compactness result obtained in Corollary \ref{comp}, it converges towards $u$ which is a solution of (\ref{eqnonh}). Since
$u\geq m\phi$,  it is non trivial.

\bigskip
We now consider the case where $\beta $ changes sign. We begin to
construct a subsolution. 

Let $\Omega^+=\{x\in\Omega\ \mbox{such that}\ \beta(x)>0\}$ which by
hypothesis is not empty. The previous case ensures that there exists a
non negative solution, not identically zero, denoted $\tilde u_0$
such that 
$$\left\{\begin{array}{lc}
G(x,\nabla \tilde u_0, D^2 \tilde u_0)  +(c+\lambda) (\tilde u_0)^{\alpha+1} 
= -\beta^+ (\tilde u_0)^{q+1},&\mbox{in}\ \Omega^+
\\ \tilde u_0=0,& \mbox{on
}\partial\Omega^+.\\
\end{array}
\right.
$$ 
We shall denote  by $u_0$ the extension:
$$u_0= \left\{ \begin{array}{cc}
\tilde u_0&{\rm in }\ \Omega^+\\
0&{\rm in } \ \Omega\setminus\Omega^+.
\end{array}\right.
$$
It is immediate to see  that $u_0$ is a  nonnegative viscosity subsolution of
$$\left\{\begin{array}{lc}
G(x, \nabla  u_0, D^2 u_0) +(c+\lambda) (u_0)^{\alpha+1} +\beta(u_0)^{q+1}
\geq  0,&\mbox{in}\ \Omega
\\ u_0=0,& \mbox{on
}\partial\Omega.\\
\end{array}
\right.
$$ 
Let $u_n$ be defined in a recursive way as 

$$\left\{\begin{array}{rc}
G(x, \nabla  u_n, D^2 u_n)+(c+\lambda-2|c+\lambda|_\infty|)
|u_n|^{\alpha} u_n-\beta^- |u_{n}|^q u_n = & \\
-2|\lambda+c|_\infty |u_{n-1}|^\alpha u_{n-1}-\beta^+
|u_{n-1}| ^q u_{n-1}&\mbox{in}\ \Omega
\\ u_n=0& \mbox{on
}\partial\Omega.\\
\end{array}
\right.
$$ 
Again, this sequence is well defined by Proposition \ref{prop0}.
We have 
$u_n\geq u_0$, and even more precisely that 
$u_n\geq u_{n-1}$. 

We claim that $u_n$ is bounded. Suppose by contradiction that  $|u_n|_\infty \rightarrow
+\infty$, then defining $w_n = {u_n\over |u_n|_\infty}$,
 one easily obtains,
as in the previous case, that $w_n$ converges, up to  
a subsequence, to some
function $w$ which satisfies for some $k\leq 1$:
$$ G(x, \nabla  w, D^2 w)+ (c+\lambda-2|c+\lambda|_\infty(1-k^{1+\alpha}))
w^{1+\alpha} =0$$
and is zero on the boundary.  One gets a contradiction with the maximum
principle and then the sequence $(u_n)$ is bounded. Extracting from it a
subsequence, using Corollary \ref{comp} and passing to the limit one gets
that 
$u$ is  a solution. 
Since $u\geq u_0$ which is not identically zero, we get the result.
 
\section{Maximum and comparison principles}

\begin{theorem}\label{maxp2}

Suppose that $F$ satisfies (H1), (H2), (H4),  that $b$ and $c$ are
continuous and
$b$ satisfies (H5).  Suppose that $h$ 
is a continuous function such that $h(x,.)$ is non  increasing, 
$h(x,0)=0$. Suppose that  $\tau<\bar\lambda$ and 
$u$ is a  viscosity sub solution of 

$$G(x,  \nabla u, D^2 u)
+h(x,u)+ ( \tau+c(x)) |u|^{\alpha}u \geq
0\ 
\ {\rm in}\
\
\Omega.$$

If $u\leq 0$ on the boundary of $\Omega$, then 
$u\leq 0$ in $\Omega$. 

\end{theorem}

{\bf Remark:} Similarly it is possible to prove that if
$\tau<
\underline\lambda$ and  $v$ is a super solution of 
$$G(x,  \nabla v, D^2 v)+ (\tau+c(x)) |v|^{\alpha}v +h(x,v)\leq 0\  \
{\rm in}\
\
\Omega.$$
If  $v\geq 0$ on the boundary of $\Omega$ then $v\geq 0$ in
$\Omega$.

\begin{corollary}

Suppose that $c^\prime$ is some continuous function, $c^\prime <
c+\bar\lambda_c$. Then 
the maximum principle holds for the equation 
$$G(x,\nabla v, D^2 v)+ c^\prime |v|^\alpha v=0$$
or equivalently if $v$ is some solution of 
$$G(x,\nabla v, D^2 v)+  c^\prime |v|^\alpha v \geq 0$$
and $v\leq 0$ on the boundary, then $v\leq 0$ in $\Omega$ 
\end{corollary}

For the proof  of the corollary it is sufficient to use the previous result with $h(x,u)=(c^\prime-(c+\lambda))|u|^\alpha u $,
where
 $\lambda<\bar\lambda$ and 
$c^\prime \leq  c+\lambda$.

Before starting the proof let us remind  two results proved in \cite{BD7}:

\begin{proposition} \label{ah} Suppose that $F$ satisfies (H1) and (H2), and
that 
$b$ is bounded.

Let $u$ be  uppersemicontinuous  subsolution of 

$$\left\{
\begin{array}{lc}
F(x, \grad u,D^2u)+b(x).\grad u |\grad u|^{\alpha}\geq
 -m &
{\rm in}\ \Omega\\
u =0 & {\rm on}\ \partial\Omega
\end{array}
\right.
$$ 
for some constant $m\geq 0$. 
Then there exists $\delta>0$ and some constant $C_3$  that depends only on the structural
data such that, for $x$ satisfying $d(x)\leq \delta$, $u$ satisfies
 $$u(x)\leq C_3d(x).$$
\end{proposition}
\begin{proposition}[Hopf] \label{Hopf}Let
$v$ be a viscosity continuous super solution of 
$$F(x,\nabla v, D^2 v)+ b(x).\nabla v |\nabla v|^\alpha +c(x) |v|^{\alpha}v
\leq 0.$$
Suppose that $v$ is positive in a neighborhood of $x_o\in\partial\Omega$ and
$v(x_o)=0$ then there exist $C>0$  and $\delta>0$ such that
$$v(x)\geq C|x-x_o|$$

for $|x-x_o|\leq \delta$.
\end{proposition}
{\bf Proof of Theorem \ref{maxp2}.}

Let $\lambda\in ]\tau, \bar\lambda[ $, and let 
$v$ be a   super solution of 
$$G(x,\nabla v, D^2v) +(\lambda +c(x))v^{\alpha+1}\leq 0,$$
satisfying $v>0$ in $\Omega$, which exists by definition of $\bar\lambda$.

We assume by contradiction that $\sup u(x)>0$ in $\Omega$. 
We first want to prove that  $\sup\frac{u}{v}<+\infty$.

Let $d(x)=d(x,\partial\Omega)$ the distance from the boundary of $\Omega$.

 Clearly from Proposition
\ref{ah}  and 
\ref{Hopf}  there exists $\delta>0$ such that $u(x)\leq Cd(x)$ and
$v(x)\geq C^\prime d(x)$ for $d(x)\leq\delta$, for some constants $C$ and
$C^\prime$.  

In the interior we just use the fact that 
$v\geq C^\prime\delta>0$ in  $\Omega_\delta=\{x:\ d(x)\geq \delta\}$ and we can conclude that 
 ${u\over v}$  is
bounded in $\Omega$.

We now define $\gamma^\prime =\sup_{x\in {\Omega}}{u\over v}$ achieved on
some point $\bar y$  and
$w(x)=\gamma v(x)$, where $0<\gamma<\gamma^\prime$,  and $\gamma$ is 
sufficiently close to $\gamma^\prime$ 
in order that 
$\displaystyle{\lambda-\tau \left({\gamma^\prime\over \gamma}\right)^{1+\alpha}\over
\left({\gamma^\prime\over \gamma}\right)^{1+\alpha}-1}\geq 2|c|_\infty$.
Furthermore by definition of the supremum there exists $x\in\Omega$ such
that
$\frac{u(x)}{v(x)}\geq \gamma$.

Clearly, by homogeneity, $G(x,\grad w,D^2 w)+(c(x)+\lambda) w^{1+\alpha}\leq 0$. 

The supremum of $u-w$ is strictly positive, and it is necessarily 
achieved inside $\Omega$, -say on $\bar x$-  since on the boundary
$u-w\leq 0$.

Let us note that 

$$(u-w)(\bar x)\geq (u-w)(\bar y)=( \gamma^\prime-\gamma )
v(\bar y)$$
and 
$$(u-w)(\bar x)\leq (\gamma^\prime-\gamma )v(\bar x) $$
which implies that 
$${\gamma^\prime\over \gamma} w(\bar x)\geq u(\bar x).$$

 We consider, for $j\in \N$ and for some $q> \max ( 2, {\alpha+2\over
\alpha+1})$: 
 $$\psi_j (x,y)
= u(x)-w(y)-{j\over q} |x-y|^q.$$  
Since $\sup (u-w)>0$,  the supremum of $\psi_j$ is achieved on 
$(x_j,y_j)\in\Omega^2$. It is classical that 
for $j$ large enough,
$\psi_j$ achieves its positive maximum on some couple
$(x_j, y_j)\in
\Omega^2$ such that 

1) $x_j\neq y_j$ for $j$ large enough,  

2) $(x_j, y_j)\rightarrow (\bar x,\bar x)$ which is a maximum point for ${u - w}$ and $\bar x$ is an interior point

3) $j|x_j-y_j|^q\rightarrow 0$,

4) there exist $X_j$ and $Y_j$ in $S$ such that 
$$\left(j|x_j-y_j|^{q-2} (x_j-y_j), X_j\right)\in J^{2,+} u(x_j)$$
and 
$$\left(j|x_j-y_j|^{q-2} (x_j-y_j), -Y_j\right)\in J^{2,-} w(y_j).$$
Furthermore
$$\left(\begin{array}{cc}
X_j&0\\
0&Y_j
\end{array}\right)\leq j\left(\begin{array}{cc}
D_j&-D_j\\
-D_j&D_j
\end{array}\right)\leq 2^{q-2} j q(q-1) |x_j-y_j|^{q-2}
\left(\begin{array}{cc}
I&-I\\
-I&I\end{array}\right)
$$ with 
$$D_j = 2^{q-3} q|x_j-y_j|^{q-2} (I+{(q-2)\over |x_j-y_j|^2}
(x_j-y_j)\otimes (x_j-y_j)).$$ 
The proof of these facts proceeds similarly to the one given in  \cite{BD7}.

\bigskip

Condition (H4) implies that

$$F(x_j, j(x_j-y_j)|x_j-y_j|^{q-2} , {X_j})-F(y_j, 
j(x_j-y_j)|x_j-y_j|^{q-2}
, {-Y_j})\leq
\omega (j|x_j-y_j|^q).$$

 Then, using the above inequality, the properties of the 
sequence $(x_j,y_j)$, the condition on $b$ -with $C_b$ below being either
the H\"older constant or $0$-, and the homogeneity condition (H1) one
obtains 

\begin{eqnarray*} 
 -h(x_j, u(x_j))-(\tau+c(x_j))u (x_j)^{1+\alpha} 
 & \leq &
G(x_j,j(x_j-y_j)|x_j-y_j|^{q-2} , X_j)\\ 
 &\leq &G(y_j, 
j(x_j-y_j)|x_j-y_j|^{q-2}, -Y_j)\\
& + & {\omega (j|x_j-y_j|^q)} + C_bj^{1+\alpha}|x_j-y_j|^{q(1+\alpha)}+o(1)\\
&\leq &-(\lambda+c(y_j)) w(y_j)^{1+\alpha}+o(1).
\end{eqnarray*}

By passing to the limit when $j$ goes to infinity, since $c$ is continuous one gets

$$-(\tau+c(\bar x)) u(\bar x)^{1+\alpha}- h(\bar x, u(\bar x))\leq-(\lambda
+c(\bar x))w(\bar x)^{1+\alpha}$$
Suppose first that $c(\bar x)+\lambda >0$ then
 using the inequality on
$w$ and the fact that $-h(\bar x, u(\bar x))\geq 0$, one gets 

$$-(\tau+c(\bar x)) u^{1+\alpha} (\bar x) \leq -(\lambda+ c(\bar x))
w^{1+\alpha}(x)$$
 
and then 

$$-(\tau+c(\bar x)) u(\bar x)^{1+\alpha} \leq -(\lambda +c(\bar x))
\left({\gamma\over \gamma^\prime} u(\bar x)\right)^{1+\alpha}, $$
which is a contradiction with the assumption on  $\gamma$ and 
$\gamma^\prime$. 

If  $(\lambda+c(\bar x))=0$  then  $\tau< \lambda$
implies that 
$-(\tau+c(\bar x))>0$ and then 
$$0<-(\tau+c(\bar x)) u^{1+\alpha} (\bar x)\leq -(\lambda+c(\bar x))\left({\gamma\over \gamma^\prime} u(\bar x)\right)^{1+\alpha}=
0$$
a contradiction. 

Suppose finally that  $\lambda+c(\bar x)<0$
then, using 
$$w(\bar x) \leq u(\bar x)+ ({\gamma\over \gamma^\prime }-1)u(\bar y)
\leq  u(\bar x)
$$
we get

$$-(c(\bar x) +\tau)u(\bar x)^{1+\alpha} \leq -(c(\bar x) +\lambda) u(\bar
x)^{1+\alpha}. $$
This implies implies that 
$$(\lambda-\tau) u(\bar x)^{1+\alpha} \leq 0$$
 a contradiction. 
This ends the proof.

\bigskip

We now prove a comparison result.

\begin{theorem}\label{complambda}
Suppose that $F$ satisfies (H1), (H2), and (H4),  that $b$ and $c$ are
continuous and
bounded and
$b$ satisfies
$(H5)$.  Suppose that $\tau< \bar\lambda$, $f_1\leq 0$, $f_1$ is
upper semi-continuous and
$f_2$ is lower semi-continuous  with $f_1\leq  f_2$. Suppose that $h$ is such that,
for all $x\in \Omega$,
$t\mapsto {-h(x,t)\over t^{\alpha+1}}$ is non decreasing on $\R^+$.

\noindent 
 Suppose that there exist $\sigma$, and $v$ non-negative, with 
$$G(x, \grad v,D^2v)+(c(x)+\tau) v^{1+\alpha}+h(x,v)\leq f_1\ \mbox{in}\quad \Omega$$
$$G(x,  \grad \sigma,D^2\sigma) 
+(c(x)+\tau)|\sigma|^{\alpha}\sigma +h(x,\sigma) 
\geq f_2\  \mbox{in}\quad \Omega$$
$$ \sigma\leq v\  \mbox{on}\quad \partial\Omega$$
Then $\sigma \leq v$ in
$\Omega$ in each of these two cases

\noindent 1) If $v>0$ on $\overline{ \Omega}$ and either $f_1<0$ in 
$\Omega$, 
 or  $f_2(\bar x)>0$ on every point $\bar x$ such that  $f_1(\bar x)=0$,

\noindent  2) If $v>0$ in $\Omega$, $f_1<0$ and $f_1<f_2 $ on
$\overline\Omega$ 
\end{theorem}

\begin{remark} Of course a similar comparison principle can be proved  for  $\tau< \underline\lambda$ and non positive solutions.
\end{remark}

\begin{corollary} \label{cc}

Suppose that $c+\lambda >0$
and $\lambda<\overline\lambda$. 
Let $u$ and $v$ be two solutions of 
$$G(x, \nabla u, D^2 u)+(c+\lambda)|u|^\alpha u = 0\quad  \mbox{in}\ \Omega.$$
1) If $u\geq v$ on $\partial \Omega$

Then $u\geq v$ in $\Omega$. 

2) If $u> v> 0$ on $\partial \Omega$ then $u> v$ in $\Omega$. 

\noindent In particular this implies that if $g>0$ the solution of
$$\left\{\begin{array}{lc}
G(x, \nabla u, D^2 u)+(c+\lambda)|u|^\alpha u = 0& \mbox{in}\ \Omega\\
u=g &\ \mbox{on} \ \partial \Omega
\end{array}
\right.$$
is unique.
\end{corollary}
{\bf Proof of Corollary \ref{cc}}
Let $u_\epsilon = u-\epsilon$ and $v_\epsilon = {v\over 1+\gamma
\epsilon}$  with $\epsilon$ and $\gamma$ chosen conveniently in order that 
$u_\epsilon \geq v_\epsilon$ on $\partial\Omega$. They satisfy
$$G(x,\nabla u_\epsilon,D^2 u_\epsilon) + (\lambda+c)(u_\epsilon)^{1+\alpha} <0$$
and 
$$G(x,\nabla v_\epsilon,D^2 v_\epsilon) + (\lambda+c)(v_\epsilon)^{1+\alpha}=0. $$
Hence  we are in the hypothesis 
of  the comparison theorem and $u_\epsilon\geq v_\epsilon$ in $\Omega$. Passing to the limit one gets the result.

We now treat the strict comparison principle.
Since $u> v>0$ there exists $\epsilon$ such that 
$u\geq (1+\epsilon ) v$ on the boundary. Since $v(1+\epsilon)$ is still a
solution by homogeneity, one gets by the first part of the corollary that 
$u\geq (1+\epsilon) v$ in $\Omega$ and then $u> v$. 

\bigskip

\noindent{\bf Proof of Theorem \ref{complambda}.}

We act as in the proof of Theorem 3.6 in \cite{BD1}.

1) We assume first that $v>0$ on $\overline{\Omega}$. 

Suppose by contradiction that $\sigma > v$ somewhere in $\Omega$. The
supremum of the function $\displaystyle{\sigma\over v}$  on
$\partial
\Omega$ is  less  than $1$ since $\sigma\leq v$ on $\partial \Omega$ and
$v>0$ on $\partial \Omega$, then its supremum is achieved inside
$\Omega$. Let
$\bar x$ be a point such that 
$$1<{\sigma(\bar x)\over v(\bar x)}= \sup_{x\in \overline{\Omega}}
{\sigma(x)\over v(x)}.$$

We define 

$$\psi_j(x,y) = {\sigma(x)\over v(y)}-{j\over qv(y)} |x-y|^q.$$
For $j$ large enough, this function achieves its maximum which is greater
than 1, on some couple $(x_j, y_j)\in
\Omega^2$. It is easy to see that this sequence converges to
$(\bar x,
\bar x)$,  a maximum point for ${\sigma\over v}$. We prove as in
\cite{BD7}  that  $x_j$, $ y_j$ can be chosen such that $x_j\neq y_j$ for
$j$ large enough.

Moreover there exist $X_j$ and $Y_j$ such that 

$$\left(j|x_j-y_j|^{q-2} (x_j-y_j), {X_j\over v(y_j)}\right)\in J^{2,+}
\sigma(x_j)$$

and

$$\left(j|x_j-y_j|^{q-2} (x_j-y_j){v(y_j)\over \beta_j}, {-Y_j\over \beta_j}\right)\in J^{2,-} v(y_j)$$
 where $\beta_j = \sigma(x_j)-{j\over q} |x_j-y_j|^q$

and 
$$F(x_j, j|x_j-y_j|^{q-2} (x_j-y_j), X_j)- F(y_j, j|x_j-y_j|^{q-2}
(x_j-y_j), {-Y_j})\leq \omega (v(y_j)j|x_j-y_j|^q).$$

We can use the fact that $\sigma$ and $v$ are respectively sub and super solution to obtain:

\begin{eqnarray*}
f_2(x_j)-(\tau +c(x_j)) \sigma
(x_j)^{1+\alpha}&-&h(x_j,\sigma (x_j))\leq G(x_j,
j|x_j-y_j|^{q-2} (x_j-y_j), {X_j\over v(y_j)})\\
 &\leq & {\beta_j^{1+\alpha}\over v(y_j)^{1+\alpha}} 
\Bigl\{F(y_j, j|x_j-y_j|^{q-2}  (x_j-y_j){v(y_j)\over \beta_j}, {-Y_j\over
\beta_j})\\
&&+\omega(jv(y_j)|x_j-y_j|^q)\Bigr\}+\\
&& + b(x_j).
j^{1+\alpha}|x_j-y_j|^{(q-1)(1+\alpha)-1} (x_j-y_j)\\ 
&\leq&
{\beta_j^{1+\alpha}\over v(y_j)^{1+\alpha}} 
G(y_j, j|x_j-y_j|^{q-2}  (x_j-y_j){v(y_j)\over \beta_j}, {-Y_j\over
\beta_j})\\
&& +{\omega(v(y_j)j|x_j-y_j|^q)\over
v(y_j)^{1+\alpha}} +C(j|x_j-y_j|^q)^{1+\alpha}\\
 &\leq &
(-\tau-c(y_j)) \beta_j^{1+\alpha} -{\beta_j^{1+\alpha}
\over v(y_j)^{1+\alpha}}h(y_j,v(y_j))+\\
&&  {\beta_j^{1+\alpha}\over
v(y_j)^{1+\alpha}} f_1(y_j)+o(1).
\end{eqnarray*}
Passing to the limit, since $c$ is continuous, we get:

$$f_2(\bar x)\leq  
\left(\frac{\sigma(\bar x)}{v(\bar x)}\right)^{\alpha+1} f_1(\bar x)
+h(\bar x,\sigma (\bar x))-\left({\sigma (\bar x)\over v(\bar x)}\right)^{1+\alpha}
h(\bar x,v(\bar x))\leq
\left(\frac{\sigma(\bar x)}{v(\bar x)}\right)^{\alpha+1} f_1(\bar x)$$
By the hypothesis on $h$.

$$f_1(\bar x)\left[1-\left(\frac{\sigma(\bar x)}{v(\bar x)}\right)^{\alpha+1}\right] \leq f_1(\bar x)-f_2(\bar x)\leq 0.$$
This contradicts the hypothesis on $f$ and $g$.

\noindent 2) Suppose that  $v>0$ and let $m>0$ be such that $f_1<f_2-m$ in $\bar\Omega$.

Using continuity there exists $\varepsilon_o>0$ such that for any $\varepsilon\leq \varepsilon_o$ :
$$ |c+\lambda|_\infty((v+\varepsilon)^{\alpha+1}-v^{\alpha+1})+|h(x,v+\varepsilon)-h(x,v)|<\frac{m}{2}.
$$

Then $w:=v+\varepsilon$ is a solution of
$$G(x,\grad w,D^2 w)+(c(x)+\lambda)w^{1+\alpha} +h(x,w)\leq f_1(x)+\frac{m}{2}< f_2(x)$$
furthermore $w\geq \sigma$ on $\partial \Omega$. And we can conclude using 
the first part that $v+\epsilon\geq \sigma$. Letting $\varepsilon$ go to zero we obtain the required result.

{\bf Remark:} The condition on the increasing behavior of $h$ in Theorem \ref{complambda} 
is somehow optimal in the sense that it is 
possible to construct a counter example when $h(x,t)\over
t^{\alpha+1}$  is non increasing :

\begin{proposition}
Let
$q<\alpha$ and
$\beta >0$, $\lambda >0$. There exists $\epsilon_2>0$ and $M>0$ such
 that for any continuous function $f$ satisfying  $0>f(x)>-\epsilon_2$ there exist two solutions $u$ and
$v$ of
$$\left\{
 \begin{array}{lc}
 G(x,\grad u,D^2 u)-\beta u^{1+q}+\lambda u^{1+\alpha}=f(x) & \mbox{in }\
\Omega\\ u=M & \mbox{on }\ \partial\Omega
\end{array}
\right.
$$
with $u\leq M\leq v$, $u\not\equiv M$ and $v\not\equiv M$.
\end{proposition}
{\bf Proof.}

Let $k(t) = -\beta t^{q+1}+\lambda
t^{1+\alpha}$, defined on $\R^+$. Let
$M>0$ and $M^\prime$ be defined as 
$$k^\prime (M)=0$$
$$k(M^\prime)=0.$$
On has $M< M^\prime$. 
Let $\epsilon_2>0$ such that 
$$k(M)=-\epsilon_2.$$
Let  $m(x)\in (M,M^\prime)$ such that $k(m(x))=f(x)>-\epsilon_2$.

We define first a sequence $(u_n)$ for $n\geq 1$:
$$\left\{\begin{array}{lc}
G(x,\nabla u_n , D^2 u_n)-\beta u_n^{1+q} = f -\lambda
u_{n-1}^{1+\alpha} & \mbox{in}\ \Omega\\
u_n = M\ & {\rm on} \ \partial \Omega;
\end{array}
\right.$$
initializing with $u_0=0$. One easily has $u_n\geq 0$ and $u_n\leq M$.
Indeed 
if $u_{n-1} \leq M$
\begin{eqnarray*}
G(x,\nabla u_n, D^2 u_n)-\beta u_n^{1+q} &\geq&  k(m(x))-\lambda 
M^{1+\alpha}\\
&  = &k(m(x))-\beta  M^{1+q}-k(M)\\
&> & f(x)-\beta M^{1+q}-\epsilon_2 >-\beta M^{1+q}.
\end{eqnarray*}
On the other hand, by definition, $M$ is a super solution and then
$$u_n\leq M$$
since on the boundary 
$u_n= M$.

In the same manner one can check that 
$u_n$ is increasing , and by passing to the limit one gets a solution. 
This solution is between $0$ and $M$ and cannot be equal to $M$.

We define next a sequence of solutions of 
$$\left\{\begin{array}{lc}
G(x,\nabla v_n, D^2 v_n)-\beta v_n^{1+q} = k(m(x)) -\lambda
v_{n-1}^{1+\alpha} & \mbox{in}\ \Omega\\
v_n = M\ & {\rm on} \ \partial \Omega
\end{array}
\right.
$$
initializing with 
$v_0 = M^\prime$. 

Let us prove that 
$v_n\geq M$:

\begin{eqnarray*}
G(x, \nabla v_1, D^2 v_1)-\beta v_1^{1+q}&=
& k(m(x))-\lambda(M^\prime)^{1+\alpha}\\
& =& -\beta m(x)^{1+q}+\lambda
m(x)^{1+\alpha}-\lambda(M^\prime)^{1+\alpha}
\\ 
&\leq & -\beta m(x)^{1+q}\leq -\beta M^{1+q}.
\end{eqnarray*}
by the comparison principle $v_1\geq M$. In the same manner we can prove the induction step and the result holds for any $n$.

We now prove that $v_1\leq M^\prime$ :
\begin{eqnarray*}
G(x, \nabla v_1, D^2 v_1)-\beta
v_1^{1+q}&=& k(m(x))-\lambda(M^\prime)^{1+\alpha}\\
&=& k(m(x))-k(M')-\beta(M^\prime)^{1+\alpha}\geq
-\beta(M^\prime)^{1+\alpha},
\end{eqnarray*}
 then $v_1\leq M'$. In the same manner, one
can prove the induction step and that the sequence is decreasing.

The limit is a solution which is between $M$ and $M^\prime$ and cannot
be equal to $M$. 
This ends the proof.

 \end{document}